\documentclass[12pt,a4paper,reqno]{amsart}

\usepackage{amsmath}
\usepackage{amssymb}
\usepackage{amsfonts}
\usepackage{xspace}

\advance\textheight by \topskip
\newtheorem{theorem}{Theorem}

\def\bb#1{[\![#1]\!]}


\title[Inversions and a parameter of Knuth]
{Combinatorics of geometrically distributed random variables:\\
Inversions and a parameter of Knuth}

\author{Helmut Prodinger}

\address{ Helmut Prodinger,
Centre for Applicable Analysis and Number Theory,
 Department of Mathematics,
University of the Witwatersrand, P.~O. Wits, 
2050 Johannesburg, South Africa, email:
{\tt helmut@gauss.cam.wits.ac.za}.
}

\date{March 17, 1999}

\begin{document}

\begin{abstract}
For words of length $n$, generated by independent geometric random
variables, we consider the mean and variance
of the number of inversions and of a parameter of Knuth from
permutation in situ. In this way, 
$q$--analogues for these parameters from the usual permutation model
are obtained.
\end{abstract}

\maketitle

\section{Introduction}
\label{sec:intro}
Let $X$ denote a geometrically distributed random variable, i.~e.
$\mathbb{P}\{X=k\}=pq^{k-1}$ for $k\in\mathbb{N}$ and $q=1-p$. The
combinatorics of $n$ geometrically distributed independent random variables
$X_1,\ldots,X_n$ has attracted recent interest, especially because of
applications in computer science. We mention just two areas, the skip list
\cite{Devroye92, PaMuPo92, Pugh90, KiPr94, Prodinger93b, KiMaPr95}
 and probabilistic counting
\cite{FlMa85, KiPr90, KiPr93, KiPrSz96}.

In \cite{Prodinger96c} the number of left-to-right maxima was investigated
for words $x_1\ldots x_n$, where the letters $x_i$ are independently
generated according to the geometric distribution. In \cite{KnPr00} the
study of left-to-right maxima was continued, but now the parameters studied
were the mean value and mean position of the $r$-th maximum.

In \cite{GrKnPr99} runs of consecutive equal letters in a string
of $n$ geometrically distributed independent random letters
were studied.

In the present paper we deal with the number of inversions.
This parameter is well understood in the context of 
permutations, see e.~g. \cite{Knuth98}.
An inversion in a word $x_1\dots x_n$ is a pair
$1\le i<j\le n$ such that $x_i>x_j$. In section 
\ref{sec:inv}
we compute average and variance of this parameter.
Interestingly, if we perform the limit $q\to1$ in these answers,
we get exactly the same formul{\ae} as in the model of 
permutations.

Another parameter related to pairs of indices in a permutation
is the parameter $a$ that was studied by Knuth in the
context of an algorithm to permute a file in situ \cite{Knuth72}, 
compare also \cite{KiPrTi87, SeFl96}. 
This parameter is defined as
\begin{equation*}
a=\big|\{(i,j)\mid 1\le i<j\le n,\ x_i=\min\{
x_i,x_{i+1},\dots,x_j
\}\}\big|.
\end{equation*}

In this more complicated example, surprisingly, the limiting case
$q\to1$ again gives exactly the
formul{\ae} from the model of 
permutations (see Section \ref{sec_knuth}).

Thus the examples treated in this paper can be interpreted as
$q$--analogues of the two parameters.

\section{The number of inversions}\label{sec:inv}

The probability that a random word of length $n$, produced
by indepent geometric random variables, has $k$
inversions, is given
as the coefficient of $v^k$ in 

\begin{align*}
f(v)=\Big(\frac pq\Big)^{n}
\sum_{i_1,\dots,i_n\ge1}
q^{i_1+\dots+i_n}
\prod_{1\le j<k\le n}
\Big(\bb{i_j\le i_k}+\bb{i_j>i_k}v
\Big).
\end{align*}

Here, $\bb{P}$ is a characteristic function, being 1 when condition
$P$ is satisfied and 0 otherwise. 
This is the notation of Iverson, being popularized by
\cite{GrKnPa94}.
The form of this generating function
is merely a reformulation of the definition of the number of inversions.

The expected value is obtained as $\mathbb{E}=
f^{'}(1)$, which is
\begin{align*}
\mathbb{E}&=\Big(\frac pq\Big)^{n}
\sum_{i_1,\dots,i_n\ge1}
q^{i_1+\dots+i_n}
\sum_{1\le j<k\le n}
\bb{i_j>i_k}\\
&=\binom n2\Big(\frac pq\Big)^{n}
\sum_{i_1,\dots,i_n\ge1}
q^{i_1+\dots+i_n}
\bb{i_1>i_2}\\
&=\binom n2\Big(\frac pq\Big)^{2}
\sum_{i_1,i_2\ge1}
q^{i_1+i_2}
\bb{i_1>i_2}\\
&=\binom n2\Big(\frac pq\Big)^{2}
\sum_{i_1>i_2\ge1}
q^{i_1+i_2}\\
&=\binom n2\frac pq
\sum_{i_2\ge1}
q^{2i_2}\\
&=\binom n2\frac {q}{1+q}.
\end{align*}

Now we are going to compute the second factorial moment
$\mathbb{E}^{\underline{2}}$, which is obtained by
$\mathbb{E}^{\underline{2}}=f^{''}(1)$,
 since the variance
$\mathbb{V}$ is given by $\mathbb{E}^{\underline{2}}+
\mathbb{E}-\mathbb{E}^2$;

\begin{align*}
\mathbb{E}^{\underline{2}}&=\Big(\frac pq\Big)^{n}
\sum_{i_1,\dots,i_n\ge1}
q^{i_1+\dots+i_n}
\sum_{1\le j<k\le n, 1\le l<m\le n, (j,k)\neq(l,m)}
\bb{i_j>i_k}
\bb{i_l>i_m}.
\end{align*}

There are several possibilities for $(j,k)\neq(l,m)$
to hold, yielding several contribution to $\mathbb{E}^{\underline{2}}
=\mathbb{E}_1^{\underline{2}}+\mathbb{E}_2^{\underline{2}}
+\mathbb{E}_3^{\underline{2}}+\mathbb{E}_4^{\underline{2}}$.

First, all 4 indices might be mutually different;

\begin{align*}
\mathbb{E}^{\underline{2}}_1
&=
\binom{n}{2}\binom{n-2}{2}
\Big(\frac pq\Big)^{4}
\sum_{i_1,\dots,i_4\ge1}
q^{i_1+\dots+i_4}
\bb{i_1>i_2}
\bb{i_3>i_4}\\
&=
\binom{n}{2}\binom{n-2}{2}
\Big(\frac pq\Big)^{4}
\sum_{i_1>i_2\ge1,i_3>i_4\ge1}
q^{i_1+\dots+i_4}\\
&=
\binom{n}{2}\binom{n-2}{2}
\frac{q^2}{(1+q)^2}.
\end{align*}

The second contribution stems from $j=l$, $k\neq m$:

\begin{align*}
\mathbb{E}^{\underline{2}}_2
&=2\binom{n}{3}
\Big(\frac pq\Big)^{3}
\sum_{i_1,i_2,i_3\ge1}
q^{i_1+i_2+i_3}
\bb{i_1>i_2}
\bb{i_1>i_3}\\
&=2\binom{n}{3}
\Big(\frac pq\Big)^{3}
\sum_{i_1>i_2\ge1,i_1>i_3\ge1}
q^{i_1+i_2+i_3}\\
&=2\binom{n}{3}
\frac{q(1+q^2)}{(1+q)(1+q+q^2)}.
\end{align*}

The third contribution originates from $j\neq l$, $k=m$:

\begin{align*}
\mathbb{E}^{\underline{2}}_3
&=2\binom{n}{3}
\Big(\frac pq\Big)^{3}
\sum_{i_1,i_2,i_3\ge1}
q^{i_1+i_2+i_3}
\bb{i_1>i_3}
\bb{i_2>i_3}\\
&=2\binom{n}{3}
\Big(\frac pq\Big)^{3}
\sum_{i_1>i_3\ge1,i_2>i_3\ge1}
q^{i_1+i_2+i_3}\\
&=2\binom{n}{3}
\frac{q^2}{1+q+q^2}.
\end{align*}

Finally, the two cases $j<k=l<m$ and $l<m=j<k$
can be combined by symmetry;

\begin{align*}
\mathbb{E}^{\underline{2}}_4
&=2\binom{n}{3}
\Big(\frac pq\Big)^{3}
\sum_{i_1,i_2,i_3\ge1}
q^{i_1+i_2+i_3}
\bb{i_1>i_2}
\bb{i_2>i_3}\\
&=2\binom{n}{3}
\Big(\frac pq\Big)^{3}
\sum_{i_1>i_2>i_3\ge1}
q^{i_1+i_2+i_3}\\
&=2\binom{n}{3}
\frac{q^3}{(1+q)(1+q+q^2)}.
\end{align*}

Altogether we find

\begin{align*}
\mathbb{E}^{\underline{2}}
=\frac{n(n-1)(n-2)}{12}\Big(3nq(1+q+q^2)
+3q^3+7{q}^{2}-q+4\Big)
{\frac {q}
{ (1+q)^{2} (1+q+q^2)}}.
\end{align*}

The variance is thus
\begin{align*}
\mathbb{V}=\frac{n(n-1)}{6}\frac{q}{(1+q)^2(1+q+q^2)}
\Big(2(1-q+q^2)n
-q^2+7q-1\Big).
\end{align*}

Summarizing, we obtain the following theorem.

\begin{theorem}The average and the  variance of the number
of inversions in a random word of length $n$ obtained by
independent geometric random variables with probabilities
$\mathbb{P}\{X=k\}=pq^{k-1}$, are given by

\begin{align*}
\mathbb{E}&=\frac{n(n-1)}{2}\frac{q}{1+q},\\
\mathbb{V}&=\frac{n(n-1)}{6}\frac{q}{(1+q)^2(1+q+q^2)}
\Big(2(1-q+q^2)n
-q^2+7q-1\Big).
\qed
\end{align*}

\end{theorem}

For fixed $q$ and $n\to\infty$, we find

\begin{align*}
\mathbb{E}&\sim\frac{n^2}{2}\frac{q}{1+q},\\
\mathbb{V}&\sim\frac{n^3}{3}\frac{q(1-q+q^2)}{(1+q)^2(1+q+q^2)}
.
\end{align*}

On the other hand, for 
$q=1$, our formul{\ae} turn into

\begin{align*}
\mathbb{E}&=\frac{n(n-1)}{4},\\
\mathbb{V}&=\frac{n(n-1)(2n+5)}{72}
,
\end{align*}
and these are exactly the 
formul{\ae} for the instance of permutions, compare
e.~g. \cite{Knuth98}.

\section{Knuth's parameter from permutation in situ}
\label{sec_knuth}

This time, the generating function of interest is
\begin{align*}
f(v)=\Big(\frac pq\Big)^{n}
\sum_{i_1,\dots,i_n\ge1}
q^{i_1+\dots+i_n}
\prod_{1\le j<k\le n}
\Big(\bb{i_j=\min \{i_j,\dots,i_k\}}v+\bb
{i_j\not=\min \{i_j,\dots,i_k\}}
\Big).
\end{align*}
Again, this is not really a useful generating function, but merely
a direct translation of the definition. Nevertheless we find it
appropriate in order to control the rather unwieldy expression.

As always,  the expected value is again obtained via
$\mathbb{E}=f^{'}(1)$;

\begin{align*}
\mathbb{E}&=\Big(\frac pq\Big)^{n}
\sum_{i_1,\dots,i_n\ge1}
q^{i_1+\dots+i_n}
\sum_{1\le j<k\le n}
\bb{i_j=\min \{i_j,\dots,i_k\}}\\
&=\sum_{1\le j<k\le n}
\Big(\frac pq\Big)^{k+1-j}
\sum_{i_j=\min\{i_j,\dots,i_k\}}q^{i_j+\dots+q_k}\\
&=\sum_{1\le j<k\le n}
\Big(\frac pq\Big)^{k+1-j}
\sum_{i\ge1}q^{i(k+1-j)}\frac{1}{p^{k-j}}\\
&=p\sum_{1\le j<k\le n}
\Big(\frac 1q\Big)^{k+1-j}
\sum_{i\ge1}q^{i(k+1-j)}\\
&=p\sum_{1\le j<k\le n}
\frac{1}{1-q^{k+1-j}}\\
&=p\sum_{2\le h\le n}(n+1-h)
\frac{1}{1-q^{h}}.
\end{align*}

And the second factorial moment is again obtained by
a second derivative;

\begin{multline*}
\mathbb{E}^{\underline{2}}=\Big(\frac pq\Big)^{n}
\sum_{i_1,\dots,i_n\ge1}
q^{i_1+\dots+i_n}\times\\ \times
\sum_{1\le j<k\le n, 1\le l<m\le n, (j,k)\neq(l,m)}
\bb{i_j=\min \{i_j,\dots,i_k\}}
\bb{i_l=\min \{i_l,\dots,i_m\}}.
\end{multline*}

Now there are even  more cases to be considered.
We might have disjoint intervals, overlapping intervals or
one interval being included in the other. Or, two indices might
coincide, resulting in two intervals glued together or
again one interval being included in the other with either
a common left or right endpoint. 

Assume first that 
{$1\le j<k<l<m\le n$}. The corresponding contribution
turns out to be

\begin{equation*}
p^2\sum_{1\le j<k<l<m\le n}\frac{1}{1-q^{k+1-j}}
\frac{1}{1-q^{m+1-l}}.
\end{equation*}

Observe that in general

\begin{equation*}
\sum_{1\le j<k<l<m\le n}a_{{k+1-j}}
a_{{m+1-l}}=
\sum_{2\le i,j\le n-2;\  i+j\le n}a_{i}a_{j}\binom{n+2-i-j}{2}.
\end{equation*}

The next range is given by
$1\le j<l<m<k\le n$, with a 
contribution

\begin{equation*}
p^2\sum_{1\le j<l<m<k\le n}\frac{1}{1-q^{k+1-j}}
\frac{1}{1-q^{m+1-l}}.
\end{equation*}

We have the general formula

\begin{equation*}
\sum_{1\le j<l<m<k\le n}a_{{{k+1-j}}}
a_{{m+1-l}}=\sum_{2\le i<j\le n}
a_{i}a_{j}(n+1-j)(j-i-1).
\end{equation*}

For the range
$1\le j<l<k<m\le n$ we obtain the contribution

\begin{equation*}
p^2\sum_{1\le j<l<m\le n}\frac{1}{1-q^{m+1-j}}
\frac{1}{1-q^{m+1-l}}(m-l-1).
\end{equation*}

Observe again that in general

\begin{equation*}
\sum_{1\le j<l<m\le n}a_{{{m+1-j}}}
a_{{m+1-l}}(m-l-1)=\sum_{3\le i<j\le n}
a_{i}a_{j}(n+1-j)(i-2).
\end{equation*}

Now the first range with 3 indices involved is
{$1\le j<k=l<m\le n$} with a contribution

\begin{equation*}
p^2\sum_{1\le j<k<m\le n}\frac{1}{1-q^{m+1-j}}
\frac{1}{1-q^{m+1-k}}.
\end{equation*}

Again, such a sum can be rearranged in general;

\begin{equation*}
\sum_{1\le j<k<m\le n}a_{{m+1-j}}
a_{m+1-k}=
\sum_{2\le i<j\le n}a_{i}a_{j}(n+1-j).
\end{equation*}

The next range $1\le j<l<m=k\le n$ gives a contribution

\begin{equation*}
p^2\sum_{2\le j<l<m\le n}
\frac{1}{1-q^{m+1-j}}
\frac{1}{1-q^{m+1-l}}.
\end{equation*}

Here we note also a general formula;

\begin{equation*}
\sum_{1\le j<l<m\le n}
a_{{m+1-j}}
a_{m+1-l}=
\sum_{2\le i<j\le n}a_{i}a_{j}(n+1-j).
\end{equation*}

The last range   {$1\le j=l<k<m\le n$}
gives a contribution

\begin{equation*}
p\sum_{1\le j<k<m\le n}
\frac{1}{1-q^{m+1-j}}
=p\sum_{1\le j<m\le n}
\frac{m-j-1}{1-q^{m+1-j}}.
\end{equation*}

Observe that in general

\begin{equation*}
\sum_{1\le j<m\le n}
({m-j-1})a_{m+1-j}=
\sum_{3\le i\le n}a_{i}(i-2)(n+1-i).
\end{equation*}

All these contributions come with a factor 2, because of symmetry.

Thus
\begin{align*}
\frac12\mathbb{E}^{(2)}&=
p^2\sum_{2\le i,j\le n-2;\  i+j\le n}\frac{1}{1-q^i}
\frac{1}{1-q^j}
\binom{n+2-i-j}{2}\\
&+p^2\sum_{2\le i<j\le n}\frac{1}{1-q^i}
\frac{1}{1-q^j}
(n+1-j)(j-i-1)\\
&+p^2\sum_{3\le i<j\le n}\frac{1}{1-q^i}
\frac{1}{1-q^j}
(n+1-j)(i-2)\\
&+p^2\sum_{2\le i<j\le n}\frac{1}{1-q^i}
\frac{1}{1-q^j}
(n+1-j)\\
&+p^2\sum_{2\le i<j\le n}\frac{1}{1-q^i}
\frac{1}{1-q^j}
(n+1-j)\\
&+p\sum_{3\le i\le n}\frac{1}{1-q^i}
(n+1-i)(i-2).
\end{align*}

After several tedious 
simplifications we arrive at this form;

\begin{align*}
\mathbb{E}^{(2)}&=
2p^2\sum_{1\le i<m\le n}\frac{1}{1-q^i}
\frac{1}{1-q^{m-i}}
\binom{n+2-m}{2}-2p\sum_{1\le j\le n}\frac{(n+1-j)^2}{1-q^j}
\\
&+2p^2\sum_{1\le i<j\le n}\frac{1}{1-q^i}
\frac{1}{1-q^j}
(n+1-j)(j-1)+n(n+1).
\end{align*}

However, we can still do better than that by noting that

\begin{align*}
\sum_{1\le i <m}\frac{1}{(1-q^i)(1-q^{m-i})}&=
\sum_{1\le i <m}
\left(\frac{q^i}{1-q^i}+\frac{1}{1-q^{m-i}}
\right)
\frac{1}{1-q^m}\\
&=-\frac{m-1}{1-q^m}+\frac{2}{1-q^m}
\sum_{1\le i <m}\frac{1}{1-q^i}.
\end{align*}

Thus										

\begin{align*}
\mathbb{E}^{(2)}&=
-2p^2\sum_{1\le m\le n}\frac{m-1}{1-q^m}
\binom{n+2-m}{2}-2p\sum_{1\le j\le n}\frac{(n+1-j)^2}{1-q^j}
\\
&+2p^2(n+1)\sum_{1\le i<j\le n}\frac{n+1-j}{(1-q^i)(1-q^j)}
+n(n+1).
\end{align*}

\begin{theorem}The average and the  variance of Knuth's
parameter from the permutation in situ problem for
random words of length $n$ obtained by
independent geometric random variables with probabilities
$\mathbb{P}\{X=k\}=pq^{k-1}$, are given by

\begin{align*}
\mathbb{E}&=
p\sum_{1\le i\le n}\frac{n+1-i}{1-q^i}
-n,\\
\mathbb{V}&=
-2p^2\sum_{1\le m\le n}\frac{m-1}{1-q^m}
\binom{n+2-m}{2}+2p^2\sum_{1\le i<j\le n}\frac{i(n+1-j)}{(1-q^i)(1-q^j)}
\\
&-p^2\sum_{1\le i\le n}\frac{(n+1-i)^2}{(1-q^i)^2}
+p\sum_{1\le i\le n}\frac{(n+1-i)(2i-1)}{1-q^i}.
\qed
\end{align*}

\end{theorem}
			
As a corollary, let us evaluate these quantities for fixed $q$
and $n\to\infty$.

For this purpose, we need two infinite series:
\begin{align*}
\alpha :=\sum_{i\ge1}\frac{1}{q^{-i}-1},\qquad
\beta:=\sum_{i\ge1}\frac{1}{(q^{-i}-1)^2}.
\end{align*}

Then

\begin{align*}
\mathbb{E}&=p\binom{n+1}{2}+pn\sum_{i\ge1}
\frac{1}{q^{-i}-1}-n+\mathcal{O}(1)\\
&=\frac{p}{2}n^2+\Big(\frac p2-1+p\alpha
\Big)n+\mathcal{O}(1).
\end{align*}

For the variance, the computations are a bit more
complicated.
We treat the sums separately:
\begin{align*}
\sum_{1\le m\le n}\frac{m-1}{1-q^m}
\binom{n+2-m}{2}&=\frac{n^4}{24}+\frac{n^3}{12}
-\frac{n^2}{24}+\frac{n^2}{2}
\sum_{m\ge1}\frac{m-1}{q^{-m}-1}+
\mathcal{O}(n)
\\&=\frac{n^4}{24}+\frac{n^3}{12}
-\frac{n^2}{24}+\frac{n^2}{2}\beta+
\mathcal{O}(n)
;
\end{align*}

\begin{align*}
\sum_{1\le i<j\le n}\frac{i(n+1-j)}{(1-q^i)(1-q^j)}=\frac{n^4}{24}+\frac{n^3}{12}-
\frac{n^2}{24}+\frac{n^2}{2}(\alpha+\beta)
+\mathcal{O}(n);
\end{align*}

\begin{align*}
\sum_{1\le i\le n}\frac{(n+1-i)^2}{(1-q^i)^2}
&=\frac{n^3}{3}+\frac{n^2}{2}
+2\sum_{1\le i\le n}\frac{(n+1-i)^2}{q^{-i}-1}
+\sum_{1\le i\le n}\frac{(n+1-i)^2}{(q^{-i}-1)^2}
+\mathcal{O}(n)
\\&=\frac{n^3}{3}+\frac{n^2}{2}
+2n^2\alpha
+n^2\beta+
\mathcal{O}(n)
;
\end{align*}

\begin{align*}
\sum_{1\le i\le n}\frac{(n+1-i)(2i-1)}{1-q^i}=
		\frac{n^3}{3}+\frac{n^2}{2}
+\mathcal{O}(n).
\end{align*}

Collecting we find
\begin{equation*}
\mathbb{V}=\frac{pq}{3}n^3+\Big(\frac{pq}2
-p^2(\alpha+\beta)\Big)n^2
+\mathcal{O}(n).
\end{equation*}

Now we consider the limit
$q\to1$. 

For the expectation we easily get
\begin{equation*}
\mathbb{E}=(n+1)H_n-2n.
\end{equation*}
For the variance we get
	\begin{align*}
\mathbb{V}&=
2\sum_{1\le i<j\le n}\frac{n+1-j}{j}
-\sum_{1\le i\le n}\frac{(n+1-i)^2}{i^2}
+\sum_{1\le i\le n}\frac{(n+1-i)(2i-1)}{i}\\
&=\sum_{1\le j\le n}\frac{(n+1-j)(4j-3)}{j}
-\sum_{1\le i\le n}\frac{(n+1-i)^2}{i^2}\\
&=-(n+1)H_n-(n+1)^2H_n^{(2)}+2n(n+2).
\end{align*}
Here we only used standard summations involving harmonic
numbers, as treated e.~g. in \cite{Knuth68, GrKn82}.
(Recall that the harmonic numbers of  first and second
order are defined by
\begin{equation*}
H_n=\sum_{1\le k\le n}\frac1{k},\qquad
H^{(2)}_n=\sum_{1\le k\le n}\frac1{k^2},
\end{equation*}
respectively.)

Thus, in the limiting case, 
expectation and variance are exactly the same as
in the permutation model, compare \cite{Knuth72, 
SeFl96, KiPrTi87}.

Some other $q$--analogues of harmonic numbers can be found
e.~g.  in \cite{AnUc85, AnCrSi97}.

\bibliographystyle{plain}


\end{document}